\documentclass[12pt]{amsart}
\usepackage{bbding}
\usepackage{dsfont}
\usepackage{graphicx}
\usepackage{amsthm,amscd}
\usepackage{calc}
\usepackage{hyperref}
\usepackage{mathrsfs}
\usepackage{mathrsfs}
\usepackage{CJK,fancyhdr,amscd}
\usepackage{anysize}
\usepackage{amsmath,amssymb,amsfonts}
\usepackage{epsfig}
\usepackage{latexsym}
\usepackage{hyperref}
\usepackage{indentfirst, latexsym}
\usepackage{graphics}
\usepackage{amsmath, bm}
\usepackage{amsfonts}
\usepackage{amssymb}%
\usepackage{color}
\usepackage{mathtools}
\usepackage{interval}
\intervalconfig {
	soft open fences ,
}
\providecommand{\U}[1]{\protect\rule{.1in}{.1in}}

\numberwithin{equation}{section}
\providecommand{\U}[1]{\protect\rule{.1in}{.1in}}

\newtheorem{theorem} {Theorem} [section]
\newtheorem{proposition}[theorem]{Proposition}
\newtheorem{corollary}  [theorem]     {Corollary}
\newtheorem{lemma}  [theorem]     {Lemma}
\newtheorem{example}  [theorem]     {Example}
\newtheorem{remark}  [theorem]     {Remark}
\newtheorem{definition}  [theorem]     {Definition}

\newcommand{\dmath}{\mathrm{d}\:\!}
\newcommand{\p}{\partial}
\newcommand{\bp}{\bar\partial}

\newcommand{\bthm}{\begin{theorem}}
	\newcommand{\ethm}{\end{theorem}}
\newcommand{\blem}{\begin{lemma}}
	\newcommand{\elem}{\end{lemma}}
\newcommand{\bprop}{\begin{proposition}}
	\newcommand{\eprop}{\end{proposition}}
\newcommand{\brmk}{\begin{remark}}
	\newcommand{\ermk}{\end{remark}}
\newcommand{\bcor}{\begin{corollary}}
	\newcommand{\ecor}{\end{corollary}}
\newcommand{\bdefi}{\begin{definition}}
	\newcommand{\edefi}{\end{definition}}
\newcommand{\beq}{\begin{equation}}
	\newcommand{\eeq}{\end{equation}}
\newcommand{\bex}{\begin{example}}
	\newcommand{\eex}{\end{example}}

\newcommand{\tOmega}{\widetilde{\Omega}}
\newcommand{\bOmega}{\overline{\Omega}{\mathstrut}}

\usepackage[lmargin=1in, rmargin=1in]{geometry}
\title{the parabolic quaternionic Monge-Amp\`{e}re type  equation on hyperK\"{a}hler manifolds}
\author{Jixiang Fu \and Xin Xu \and Dekai Zhang}
\address{Shanghai Center for Mathematical Sciences, Jiangwan Campus, Fudan University, Shanghai, 200438, China}
\email{majxfu@fudan.edu.cn}
\address{School of Mathematical Sciences, Fudan University, Shanghai, 200433, China}
\email{20110180015@fudan.edu.cn}
\address{ ${}^1$Department of Mathematics, Shanghai University, Shanghai, 200444, China. \qquad \qquad \linebreak ${}^2$Newtouch Center for Mathematics of Shanghai University, Shanghai, 200444, China.
}
\email{dkzhang@shu.edu.cn}

\begin{document}
\begin{abstract}
	We prove the long time existence and uniqueness of solution to a parabolic quaternionic Monge-Amp\`{e}re type  equation on a compact hyperK\"{a}hler manifold. We also show that after normalization, the solution converges smoothly to the unique solution of the Monge-Amp\`{e}re equation for $(n-1)$\nobreakdash-quaternionic psh functions.
\end{abstract}
\maketitle
\section{introduction}
\label{intro}
A hypercomplex manifold is a smooth manifold $M$ together with a triple $(I, J, K)$ of
complex structures satisfying the quaternionic relation $IJ = -JI = K.$ A hyperhermitian metric on a hypercomplex manifold $(M,I,J,K)$ is a Riemannian metric $g$ which is hermitian with respect to $I$, $J$ and $K$. 

On a hyperhermitian manifold $(M, I, J , K, g)$, let
$ \Omega = \omega_J - i\omega_K$ where $\omega_J$ and $\omega_K$ are the fundamental forms corresponding to $J$ and $K$ respectively. Then $g$ is called hyperK\"{a}hler (HK) if $d\Omega = 0$, and called hyperK\"{a}hler with torsion (HKT) if $\partial\Omega = 0$.  Throughout this paper we use $\p$ and $\bp$ to denote the complex partial differential  operator with respect to the complex structure $I$.

Analogous to the complex Calabi-Yau equation on  K\"{a}hler manifolds which solved by Yau \cite{yau1978ricci}, Alesker and Verbitsky introduced a quaternionic Calabi-Yau equation on hyperhermitian manifolds in  \cite{alesker2010quaternionic}
\begin{equation}
	\label{eqn:intro-q-calabi}
	\begin{split}
		(\Omega + \partial \partial_J u)^n &= e^f \Omega^n, \\
		\Omega + \partial \partial_J u &> 0,
	\end{split}
\end{equation}
where $f$ is a given smooth function on $M$ and $\partial_J := J^{-1} \circ \overline{\partial} \circ J $.  They conjectured that the equation is solvable on HKT manifolds with holomorphically trivial canonical bundle with respect to $I$ and further obtained the $C^0$ estimate in this setting \cite{alesker2010quaternionic}. Alesker \cite{alesker2013solvability} solved the equation on a flat hyperK\"{a}hler manifold. In \cite{ALESKER20171} Alesker and Shelukhin proved the $C^0$ estimate without any extra assumptions and the proof was later simplified by Sroka \cite{sroka2020}.
Recently Dinew and Sroka \cite{dinew2021hkt} solved the equation on a compact HK manifold.
Bedulli, Gentili and Vezzoni \cite{bedulli2023parabolic} considered the parabolic method.
 More partial results can be found in \cite{alesker2006plurisubharmonic, alesker2010quaternionic, BGV22, gentili2022quaternionic,  gentili23, sroka22, zhang22} and the conjecture remains open.

By adopting the techniques of Dinew and Sroka \cite{dinew2021hkt},  we considered the quaternionic form-type Calabi-Yau equation in \cite{fu2023monge} on compact HK manifolds, which is parallel to the complex case where the form-type Calabi-Yau equation was proposed by Fu-Wang-Wu \cite{fu2010form,fu2015form} and solved by Tosatti-Weinkove \cite{tosatti2017monge} on K\"ahler manifolds.

Specifically, let $(M, I, J, K, g, \Omega)$ be a hyperhermitian manifold of quaternionic dimension $n$,
and $g_0$ another hyperhermitian metric on $M$ with induced $(2, 0)$-form $\Omega_0$. Given a  smooth function $f$ on $M$, the quaternionic form-type Calabi-Yau equation is
\begin{equation}
	\label{eqn:func-type}
	\Omega_{u}^n = e^{f + b}\Omega^n
\end{equation}
in which $b$ is a uniquely determined constant,  and $\Omega_{u}$ is determined by
\begin{equation}
	\label{eqn:(n-1)power}
	\Omega_u^{n-1} = \Omega_0^{n-1} + \partial \partial_J u \wedge \Omega^{n-2}
\end{equation}
where $\Omega_0^{n-1} + \partial \partial_J u \wedge \Omega^{n-2}$ is strictly positive. In fact we solved the following Monge-Amp\`{e}re equation for $(n-1)$\nobreakdash-quaternionic psh functions which is equivalent to \eqref{eqn:func-type}.
\begin{equation}
\label{eqn:psh}
	\begin{aligned}
		\big(\Omega_h + & \frac{1}{n-1}((\frac{1}{2}\Delta_{I, g}u)\Omega - \partial \partial_J u)\big)^n
		= e^{f + b} \Omega^n \\
		\Omega_h &+ \frac{1}{n-1}((\frac{1}{2}\Delta_{I, g}u)\Omega - \partial \partial_J u) > 0,
	\end{aligned}
\end{equation}
where $\Omega_h$ a given strictly positive $(2, 0)$\nobreakdash-form with respect to $I$.

In a slightly different context, on locally flat compact HK manifolds,  Gentili and Zhang \cite{gentili2022fully} solved a class of fully non-linear elliptic equations including \eqref{eqn:psh}. They later extended their work to the parabolic setting in \cite{gentili2022fullyp}. 


In this article, we consider the parabolic version of \eqref{eqn:psh} on a compact hyperK\"ahler manifold
\begin{equation}
\label{eqn:flow}
	\frac{\partial}{\partial t} u = 
	\log{\frac{\big(\Omega_h + \frac{1}{n-1}((\frac{1}{2}\Delta_{I, g}u)\Omega - \partial \partial_J u)\big)^n}{\Omega^n}} -f,
\end{equation}
with $u(\cdot, 0) = u_0 \in C^{\infty} (M , \mathbb{R})$ satisfying
\begin{equation}\label{u0initial}
	\Omega_h + \frac{1}{n-1}((\frac{1}{2}\Delta_{I, g}u_0)\Omega - \partial \partial_J u_0) > 0.
	\end{equation}
  Our main result is as follows.
\begin{theorem}
	\label{thm:main}
	Let $(M, I, J, K, g, \Omega)$ be a compact hyperK\"{a}hler manifold of quaternionic dimension $n$, and $\Omega_h$ a strictly positive $(2, 0)$\nobreakdash-form with respect to $I$. Let $f$ be a smooth function on $M$. Then there exists a unique solution $u$ to \eqref{eqn:flow} on $M \times [0, \infty )$ with $u(\cdot,0)=u_0$ satisfying \eqref{u0initial}.
	And if we normalize $u$ by
	\begin{equation}
		\Tilde{u} := u - \frac{\int_M u \, \Omega^n \wedge \bOmega^n }{\int_M   \Omega^n \wedge \bOmega^n },
	\end{equation}
	then $\Tilde{u}$ converges smoothly to a function  $\Tilde{u}_{\infty}$ as $t \rightarrow \infty$, and $\Tilde{u}_{\infty}$ is the unique solution to
	 \eqref{eqn:psh} up to a constant $\Tilde{b} \in \mathbb{R}$.
\end{theorem}

This gives a parabolic solution to the original equation \eqref{eqn:psh}. There are plenty of results on  parabolic flows on compact complex manifolds, for example, \cite{Cao1985, Chu18, gill2011, flm2011, sw2008, zheng19}.

The article is organized as follows. In Section 2, we introduce some basic notations and useful lemmas. In Section \ref{C0}, we prove the $u_t$ and the $C^0$ estimate. We derive the $C^1$ estimate  in Section \ref{C1} and the complex Hessian estimate in Section \ref{pC2}. The Theorem \ref{thm:main} is proved in Section \ref{C2a}.

\section{Preliminaries}
\label{pre}
On a hyperhermitian manifold $(M, I, J , K, g)$ of quaternionic dimension $n$,  we denote by $\Lambda_{I}^{p,q}(M)$ the $(p,q)$\nobreakdash-forms with respect to $I$.
A form $\alpha \in \Lambda_{I}^{2k, 0} (M)$ is called $J$\nobreakdash-real if $J\alpha = \overline{\alpha}$, and denoted by $\alpha \in \Lambda_{I, \mathbb{R}}^{2k, 0} (M)$. In particular, we have $\Omega = \omega_J - i\omega_K$ is a $J$\nobreakdash-real $(2,0)$\nobreakdash-form.

\begin{definition}[\cite{fu2023monge}, Definition 2.2]
	A $J$\nobreakdash-real $(2,0)$\nobreakdash-form $\alpha$ is said to be
	positive (\emph{resp.} strictly positive) if
	$\alpha(X, \overline{X}J) \geq 0$ (\emph{resp.}  $\alpha(X, \overline{X}J)  > 0$)
	for any non-zero $(1,0)$\nobreakdash-vector $X$. We denote by $\Lambda_{I,\mathbb R}^{2,0}(M)_{>0}$
	all strictly positive $J$-real $(2,0)$-forms.
\end{definition}

Note that $\Omega$ is determined by $g$ and is strictly positive. Conversely any $\Omega \in \Lambda_{I, \mathbb{R}}^{2k, 0} (M)_{>0}$ induces a hyperhermitian metric by $g = \operatorname{Re}(\Omega(\cdot, \cdot J))$. Thus there is a bijection between strictly positive $J$\nobreakdash-real $(2,0)$\nobreakdash-forms and hyperhermitian metrics. 

\begin{definition}
	For $\chi \in \Lambda_{I, \mathbb{R}}^{2, 0}(M)$, define
	\begin{equation}
		\label{eqn:intro-Sm}
		S_m(\chi) = \frac{C_n^m \chi^m \wedge \Omega^{n-m}}{\Omega^n}
		\quad \text{for} \quad 0 \leq m \leq n.
	\end{equation}
\end{definition}
In particular for $u \in C^{\infty}(M, \mathbb{R})$ we have
\begin{equation}
	\label{eqn:C0-laplace}
	S_1(\partial\partial_J u ) = \frac{1}{2} \Delta_{I, g} u.
\end{equation}
For convenience we denote
\begin{equation}
	\label{equ:defOmegaTilda}
	\tOmega = \Omega_h + \frac{1}{n-1}(S_1(\partial\partial_J u)\Omega - \partial\partial_J u).
\end{equation}
It's easily checked that $\tOmega$ is a $J$\nobreakdash-real $(2,0)$\nobreakdash-form, thus one can define the corresponding hyperhermitian metric and the induced fundamental form by
\begin{equation}
	\label{def:omega-u}
	g_u = \operatorname{Re}(\tOmega(\cdot, \cdot J)), \quad
	\omega_u = g_u(\cdot I, \cdot).
\end{equation}

\begin{lemma}
	\begin{equation}
		\omega_u = \omega_h + \frac{1}{n-1}(S_1(\p \p_J u)\omega - \frac{1}{2}(i\p \bp u - i J \p \bp u)).
	\end{equation}
\end{lemma}
\begin{proof}
	It is showed in \cite[Proposition 3.2]{sroka22}  that 
	\begin{equation*}
		 \operatorname{Re}(\p \p_J u(\cdot I, \cdot J)) = \frac{1}{2} (i\p \bp u - i J \p \bp u).
	\end{equation*}
	Hence by definition
		\begin{align*}
			\omega_u &= g_u(\cdot I, \cdot) = \operatorname{Re}(\tOmega(\cdot I, \cdot J)) \\
							  &= \operatorname{Re}(\Omega_h(\cdot I, \cdot J)) + \frac{1}{n-1}(S_1(\p \p_J u)\operatorname{Re}(\Omega(\cdot I, \cdot J)) - \operatorname{Re}(\p \p_J u(\cdot I, \cdot J)))\\
							  &= \omega_h + \frac{1}{n-1}(S_1(\p \p_J u)\omega - \frac{1}{2}(i\p \bp u - i J \p \bp u)).\qedhere
		\end{align*}
\end{proof}

We also need the following lemma.
\begin{lemma}[\cite{fu2023monge}, Lemma 3.2]
	\begin{gather}
		\label{eq:S1ddju}
		S_1(\partial\partial_J u) = S_1(\tOmega) - S_1(\Omega_h),\\
		\label{eq:ddju}
		\partial\partial_J u = (n-1)\Omega_h - S_1(\Omega_h)\Omega + S_1(\tOmega)\Omega - (n-1)\tOmega.
	\end{gather}
\end{lemma}


\begin{remark}
	\label{rm:normal}
	On a hyperhermitian manifold $(M, I, J, K, g, \Omega)$ of quaternionic dimension $n$, we can find local $I$\nobreakdash-holomorphic geodesic coordinates such that $\Omega$ and another $J$\nobreakdash-real $(2,0)$\nobreakdash-form $\tOmega$ are simultaneously diagonalizable at a point $x \in M$, i.e.
	\begin{equation*}
		\Omega = \sum_{i=0}^{n-1} dz^{2i} \wedge dz^{2i+1}, \quad
		\tOmega = \sum_{i=0}^{n-1} \tOmega_{2i2i+1} dz^{2i} \wedge dz^{2i+1},
	\end{equation*}
	and the Christoffel symbol of $\nabla^O$ and first derivatives of $J$ vanish at $x$, i.e.
	\begin{equation*}
		J^l_{\bar k,i} =J^{\bar l}_{k,i} =J^{\bar l}_{k,\bar i}=J^l_{\bar k, \bar i} =0.
	\end{equation*}
	We call such local coordinates the normal coordinates around $x$.

\end{remark}
The linearized operator $\mathcal{P}$ of the flow \eqref{eqn:flow} is derived in the following lemma
\begin{lemma}
	The linearized operator $\mathcal{P}$ has the form:
	\begin{align}
		\mathcal{P}(v)=v_t-\frac{	A\wedge \p\p_J(v)}{ \tOmega^n},
	\end{align}
where $A = \frac{n}{n-1}\big(S_{n-1}(\tOmega)\Omega^{n-1} - \tOmega^{n-1}\big)$ and $v\in C^{2,1}(M\times [0, T))$.
\end{lemma}
\begin{proof}
	Let $w(s)$ be the variation of $u$ and $v=\frac{d}{ds}\Big |_{s=0}w(s)$. It is sufficient to compute the variation of $\tOmega^n=(\Omega_h + \frac{1}{n-1}(S_1(\partial\partial_J u)\Omega - \partial\partial_J u))^n$. We have
	\begin{align*}
		\delta(\tOmega^n)=& \frac{d}{ds}\Big|_{s=0}(\Omega_h + \frac{1}{n-1}(S_1(\partial\partial_J w(s))\Omega - \partial\partial_J w(s)))^n\\
		=&\frac{n}{n-1}\tOmega^{n-1}\wedge ( S_1(\p\p_J v)\Omega-\p\p_J v)
		\\
		=&\frac{n}{n-1}\tOmega^{n-1}\wedge \Omega \cdot \frac{n\Omega^{n-1}\wedge \p\p_J v }{\Omega^n}-\frac{n}{n-1}\tOmega^{n-1}\wedge\p\p_J v\\
		=&\frac{n}{n-1}S_{n-1}( \tOmega) {\Omega^{n-1}\wedge \p\p_J v }-\frac{n}{n-1}\tOmega^{n-1}\wedge\p\p_J v\\
		=&A\wedge \p\p_J v.
		\end{align*}
	
Then	$\mathcal{P}(v)=v_t-\delta (\log \frac{\tOmega^n}{\Omega^n})=v_t-\frac{	A\wedge \p\p_J(v)}{ \tOmega^n}$ as claimed.
\end{proof}

\section{\texorpdfstring{$u_t$ estimate and $C^0$}{C0} estimate}
\label{C0}

 We first prove the uniform estimate of $u_t$.
\begin{lemma}
	\label{ut-bound}
	Let $u$ be a solution to \eqref{eqn:flow} on $M \times [0, T)$. Then there exists a constant $C$ depending only on the fixed data $(I, J, K, g, \Omega, \Omega_h)$ and $f$ such that
	\begin{equation}
		\sup_{M \times [0, T)} \big|u_t\big| \leq C.
	\end{equation}
\end{lemma}
\begin{proof}
	
	One can see that	 $u_t$ satisfies 
	\begin{align}
		\mathcal{P}(u_t)=\frac{\p}{\p t}(u_t)-\frac{	A\wedge \p\p_J(u_t)}{ \tOmega^n}=0.
	\end{align}
	For any $T_0\in (0, T)$, by maximum principle, 
	\begin{align*}
		\max\limits_{M\times [0, T_0]}|u_t|\le& \max\limits_{M}|u_t(x, 0)|\\
		\le& \max\limits_M\Big|\log{\frac{\big(\Omega_h + \frac{1}{n-1}(S_1(\p \p_J u_0)\Omega - \partial \partial_J u_0)\big)^n}{\Omega^n}}\Big| +\max\limits_{M}|f|.
	\end{align*}
	Since $T_0$ is arbitrary, we have the desired estimate.
	
\end{proof}
Using the $C^0$ estimate for the elliptic equation, which has been proved by Sroka \cite{sroka22} and Fu-Xu-Zhang \cite{fu2023monge}, we have the following.

\begin{lemma}
	\label{neq:tC0}
		Let $u$ be a solution to \eqref{eqn:flow} on $M \times [0, T)$. Then there exists a uniform constant $C$
	depending only on the fixed data $(I, J, K, g, \Omega, \Omega_h)$ and $f$ such that
	\begin{equation}
		\sup\limits_{M\times[0,T)}|\tilde u| \leq
		\sup\limits_{t\in \interval[open right]{0}{T} } \big( 
		\sup\limits_{x \in  M} u(x,t) - \inf\limits_{x \in M} u(x, t) 
		 \big)
		\leq C.
	\end{equation}
	\end{lemma}
\begin{proof}
	The flow is equivalent to the following
		\begin{equation}
						\label{equation1}
						\tOmega^n = e^{u_t +f} \Omega^n.
					\end{equation}
Since $u_t$ is uniformly bounded, we can apply the $C^0$-estimate for the elliptic equation such that for any $t\in(0, T)$, 
\begin{align}
	|u(x,t)-\sup_M u(\cdot, t)|\le C, \quad \forall x\in M.
\end{align}

Since $\int\limits_M{\tilde u(\cdot, t)\,  \Omega^n\wedge\bOmega^n}=0$, there exists $x_0\in M$ such that $\tilde u(x_0, t)=0$.
Then we have	
	\begin{align*}
	|\tilde u(x,t)|=&\,|\tilde u(x,t)-\tilde u(x_0,t)|
	=| u(x,t)- u(x_0,t)|\\
	\le&\, | u(x,t)-\sup_M u(\cdot, t) |+| u(x_0,t)-\sup_M u(\cdot, t) |\\
	\le&\, 2C, \quad \forall x\in M.
\end{align*}	
Hence	the $C^0$ estimate follows.
\end{proof}

\section{\texorpdfstring{$C^1$}{C1} Estimate}
\label{C1}
 Although the gradient estimate is unnecessary for the proof of the main result,  we provide it as the gradient estimate for fully nonlinear equations has independent interest.
\begin{theorem}
	Let $u$ be a solution to \eqref{eqn:flow} on $M \times [0, T)$. Then there exists a constant $C$
	depending only on the fixed data $(I, J, K, g, \Omega, \Omega_h)$ and $f$ such that
	\begin{equation}
		\label{neq:C1}
	\sup\limits_{M\times[0,T)}	|du|_g \leq C.
	\end{equation}
\end{theorem}
\begin{proof}
	    A simple computation in local coordinates shows that
	\begin{equation*}
		n \partial u \wedge \partial_J u \wedge \Omega^{n-1} = \frac{1}{4} |d u|_g^2 \Omega^n.
	\end{equation*}
	Define
	\begin{equation*}
		\beta \coloneqq \frac{1}{4} |d u|_g^2.
	\end{equation*}
	Following \cite{blocki2009gradient}, we consider
	\begin{equation*}
		G = \log \beta - \varphi(\tilde u) ,
	\end{equation*}
	where $\varphi$ is a function to be determined and $\tilde u$ is the normalization of $u$  . For any $T_0\in (0,T)$, suppose $\max\limits_{M\times [0, T_0]}G=G(p_0, t_0)$ with $(p_0,t_0)\in M\times [0, T_0)$. We want to show $\beta(p_0,t_0)$ is uniformly bounded.
If $t_0=0$, we have the estimate. In the following, we assume $t_0>0 $.
	
	We choose the normal coordinates around $p_0$ (see Remark \ref{rm:normal}) and all the calculation is at $(p_0, t_0)$.
	\begin{equation*}
		\begin{split}
			0\le \partial_t G &= \frac{ \beta_t}{\beta} - \varphi'  \tilde u_t;\\
			\partial G &= \frac{\partial \beta}{\beta} - \varphi' \partial u = 0; \\
			\partial_J G &= \frac{\partial_J \beta}{\beta} - \varphi' \partial_J u = 0;\\
			\partial\partial_J G &= \frac{\partial \partial_J \beta}{\beta} -
			\frac{\partial \beta \wedge \partial_J \beta}{\beta^2} -
			\varphi'' \partial u \wedge \partial_J u -
			\varphi' \partial \partial_J u \\
			&= \frac{\partial \partial_J \beta}{\beta} -
			((\varphi')^2 + \varphi'') \partial u \wedge \partial_J u
			- \varphi' \partial \partial_J u.
		\end{split}
	\end{equation*}
	Then we have
	\begin{equation}
		\label{eqn:C1less0}
		\begin{split}
			0 \leq &\mathcal{P}(G)= G_t -  \frac{\partial \partial_J G \wedge A \wedge \bOmega^n}{\tOmega{\mathstrut}^n \wedge \bOmega^n} \\
			&=  \frac{\beta_t}{\beta} - \varphi' \tilde u_t - \frac{\partial \partial_J \beta \wedge A \wedge \bOmega^n}{\beta \tOmega{\mathstrut}^n \wedge \bOmega^n} +
			((\varphi')^2 + \varphi'') \frac{\partial u \wedge \partial_J u \wedge A \wedge \bOmega^n}{\tOmega{\mathstrut}^n \wedge \bOmega^n} +
			\varphi' \frac{\partial \partial_J u \wedge A \wedge \bOmega^n}{\tOmega{\mathstrut}^n \wedge \bOmega^n}.
		\end{split}
	\end{equation}

	We first deal with $\partial_t \beta$. By taking $\partial_t$ on both sides of
		$\beta \Omega^n = n \partial u \wedge\partial_J u \wedge \Omega^{n-1}$, we get
	\begin{equation}\label{bt}
		\beta_t = \sum_{j=0}^{2n-1} (u_{t,j} u_{\overline{j}} + u_j u_{t, \overline{j} }).
	\end{equation}
	
	We next compute $\partial \partial_J \beta$. 
	Taking $\partial_J$ on both sides of	$\beta \bOmega^n = n \overline{\partial} u \wedge \overline{\partial_J} u \wedge \bOmega^{n-1}$ and noticing $\partial_J \bOmega = 0$(since $\Omega$ is hyperK\"{a}hler), we have
	\begin{equation*}
		\partial_J \beta \wedge \bOmega^n = n \partial_J \overline{\partial} u \wedge \overline{\partial_J} u \wedge \bOmega^{n-1} -
		n \overline{\partial} u \wedge \partial_J \overline{\partial_J} u \wedge \bOmega^{n-1}.
	\end{equation*}
	Then taking $\partial$ on both sides, we obtain
	\begin{equation*}
		\begin{split}
			\partial \partial_J \beta \wedge \bOmega^n =
			&n \partial \partial_J \overline{\partial} u \wedge \overline{\partial_J} u \wedge \bOmega^{n-1} +
			n \partial_J \overline{\partial} u \wedge \partial \overline{\partial_J} u \wedge \bOmega^{n-1} \\
			&- n \partial \overline{\partial} u \wedge \partial_J \overline{\partial_J} u
			\wedge \bOmega^{n-1}
			+ n \overline{\partial} u \wedge \partial \partial_J \overline{\partial_J} u
			\wedge \bOmega^{n-1}.
		\end{split}
	\end{equation*}
	From the equation
	\begin{equation}
		\label{eqn:C1-main}
		\tOmega^n = e^{u_t +f} \Omega^n,
	\end{equation}
	by taking $\overline{\partial}$ on both sides we get
	\begin{equation*}
		\begin{split}
			n(\bp S_1(\p \p_J u)\wedge\Omega - \bp \p \p_J u)\wedge \tOmega^{n-1} &=
			(n-1)(\bp e^{u_t + f} \wedge \Omega^n - n \bp \Omega_h \wedge \tOmega^{n-1}).
		\end{split}
	\end{equation*}
The left hand side can be calculated as the following.
	\begin{equation*}
		\begin{split}
			&n(\bp S_1(\p \p_J u)\wedge\Omega - \bp \p \p_J u)\wedge \tOmega^{n-1} \\
			=\, &n(\bp S_1(\p \p_J u) \wedge \Omega^n \cdot \frac{\Omega \wedge \tOmega^{n-1}}{\Omega^n} - \bp \p \p_J u \wedge \tOmega^{n-1}) \\
			=\, &n(\bp \big(\frac{ \p \p_J u \wedge \Omega^{n-1}}{\Omega^n} \cdot \Omega^n\big) \cdot
			S_{n-1}(\tOmega) - \bp \p \p_J u \wedge \tOmega^{n-1} )\\
			=\, &(S_{n-1}(\tOmega)\Omega^{n-1} - \tOmega^{n-1})\wedge n \bp \p \p_J u\\
			=\, & (n-1)A \wedge  \bp \p \p_J u.
		\end{split}
	\end{equation*}
	Hence we obtain
	\begin{equation*}
		A \wedge n \overline{\partial} \partial \partial_J u =
		-n^2 \tOmega^{n-1} \wedge \overline{\partial} \Omega_h + n\overline{\partial} e^{u_t + f} \wedge \Omega^n.
	\end{equation*}
By taking $\overline{\partial_J}$ on both sides of \eqref{eqn:C1-main}, we obtain
	\begin{equation*}
		A \wedge n \overline{\partial_J} \partial \partial_J u =
		-n^2 \tOmega^{n-1} \wedge \overline{\partial_J} \Omega_h + n\overline{\partial_J} e^{u_t +f} \wedge \Omega^n .
	\end{equation*}
	Thus for the third term of \eqref{eqn:C1less0}, we have
	\begin{equation}
		\label{eqn:C1-1}
		\partial \partial_J \beta \wedge A \wedge \bOmega^n =
		I_1 + I_2 +
		n \partial_J \overline{\partial} u \wedge \partial \overline{\partial_J} u
		\wedge \bOmega^{n-1} \wedge A - n \partial \overline{\partial} u
		\wedge \partial_J \overline{\partial_J} u \wedge \bOmega^{n-1} \wedge A
	\end{equation}
	where
	\begin{equation*}
		\begin{split}
			I_1 &= (-n^2 \tOmega^{n-1} \wedge \overline{\partial} \Omega_h + n\overline{\partial} e^{u_t + f} \wedge \Omega^n)
			\wedge \overline{\partial_J} u \wedge \bOmega^{n-1}, \\
			I_2 &= (n^2 \tOmega^{n-1} \wedge \overline{\partial_J} \Omega_h - n\overline{\partial_J} e^{u_t +f} \wedge \Omega^n)
			\wedge \overline{\partial} u \wedge \bOmega^{n-1}.
		\end{split}
	\end{equation*}
	
	By direct computation,
	\begin{equation*}
		\begin{split}
			\partial_J \overline{\partial} u &= \sum u_{\overline{ji}} J^{-1} d\overline{z^i} \wedge d\overline{z^j}, \\
			\partial \overline{\partial_J} u &= \sum u_{ij} dz^j \wedge J^{-1} dz^i, \\
			\partial \overline{\partial} u &= \sum u_{i\overline{j}} dz^i \wedge d\overline{z^j}, \\
			\partial_J \overline{\partial_J} u &= \sum u_{i\overline{j}} J^{-1} d\overline{z^j} \wedge J^{-1} dz^i,
		\end{split}
	\end{equation*}
	 the third term of \eqref{eqn:C1-1} becomes
	\begin{equation}
		\label{eqn:C1-1-3}
		n \partial_J \overline{\partial} u \wedge \partial \overline{\partial_J} u
		\wedge \bOmega^{n-1} \wedge A =
		\frac{1}{n-1} \sum_{k=0}^{n-1} \sum_{j=0}^{2n-1} (\sum_{i\neq k}\frac{1}{\tOmega_{2i2i+1}} )(|u_{2kj}|^2 + |u_{2k+1j}|^2)
		\tOmega^n \wedge \bOmega^n;
	\end{equation}
	and the forth term
	\begin{equation}
		\label{eqn:C1-1-4}
		-n \partial \overline{\partial} u \wedge \partial_J \overline{\partial_J} u
		\wedge \bOmega^{n-1} \wedge A =
		\frac{1}{n-1}\sum_{k=0}^{n-1} \sum_{j=0}^{2n-1} (\sum_{i\neq k}\frac{1}{\tOmega_{2i2i+1}} )
		(|u_{2k\overline{j}}|^2 + |u_{2k+1\overline{j}}|^2)
		\tOmega^n \wedge \bOmega^n.
	\end{equation}
	For $I_1$ and $I_2$ we have
	\begin{equation}
		\label{eqn:C1-1-1}
		\begin{split}
			 I_1 &= -n^2 \tOmega^{n-1} \wedge \overline{\partial} \Omega_h \wedge \overline{\partial_J} u \wedge \bOmega^{n-1}  - n\overline{\partial_J} u \wedge  \overline{\partial} e^{u_t + f} \wedge \Omega^n \wedge \bOmega^{n-1} \\
			&= -  \sum_{i = 0}^{n-1} \sum_{j = 0}^{2n-1}
			\frac{(\Omega_h)_{2i2i+1, \overline{j}} u_{j}  }{\tOmega_{2i2i+1}}\tOmega^n \wedge \bOmega^n +
			\sum_{j=0}^{2n-1} u_j (u_t + f)_{\overline{j}} \tOmega^n \wedge \bOmega^n
		\end{split}
	\end{equation}
	and
	\begin{equation}
		\label{eqn:C1-1-2}
		\begin{split}
			 I_2 &= n \tOmega^{n-1} \wedge \overline{\partial_J} \Omega_h \wedge \overline{\partial} u \wedge \bOmega^{n-1}  + \overline{\partial} u \wedge  \overline{\partial_J} e^{u_t + f} \wedge \Omega^n \wedge \bOmega^{n-1} \\
			&= -  \sum_{i = 0}^{n-1} \sum_{j = 0}^{2n-1}
			\frac{(\bOmega_h)_{2i2i+1, j} u_{\overline{j}}  }{\tOmega_{2i2i+1}}\tOmega^n \wedge \bOmega^n +
			\sum_{j=0}^{2n-1} u_{\overline{j}}(u_t +f)_j\tOmega^n \wedge \bOmega^n.
		\end{split}
	\end{equation}
	Combining \eqref{eqn:C1-1-3}, \eqref{eqn:C1-1-4}, \eqref{eqn:C1-1-1}, \eqref{eqn:C1-1-2},  we obtain estimate of \eqref{eqn:C1-1}
	\begin{equation}
		\label{eqn:C1-1-local}
		\begin{split}
			\frac{\p \p_J \beta \wedge A \wedge \bOmega^n}{\beta \tOmega{\mathstrut}^n \wedge \bOmega^n}
			&= - \frac{1}{\beta} \sum_{i = 0}^{n-1} \sum_{j = 0}^{2n-1}
			\frac{(\Omega_h)_{2i2i+1, \overline{j}} u_{j} + (\bOmega_h)_{2i2i+1, j} u_{\overline{j}} }
			{\tOmega_{2i2i+1}}\\
			&+ \frac{1}{\beta} \sum_{j=0}^{2n-1}
			\big(u_j(u_t + f)_{\overline{j}} + u_{\overline{j}}(u_t + f)_{j} \big) \\
			&+ \frac{1}{(n-1)\beta} \sum_{k=0}^{n-1} \sum_{j=0}^{2n-1} \sum_{i\neq k}
			\frac{|u_{2kj}|^2 + |u_{2k+1j}|^2 + |u_{2k\overline{j}}|^2 + |u_{2k+1\overline{j}}|^2}{\tOmega_{2i2i+1}} .
		\end{split}
	\end{equation}
	Again by direct computation, the forth term of \eqref{eqn:C1less0} is
	\begin{equation}
		\label{eqn:C1-2}
		\partial u \wedge \partial_J u \wedge A \wedge \bOmega^n =
		\frac{1}{n-1}\sum_{i=0}^{n-1}  (\sum_{k\neq i}\frac{1}{\tOmega_{2k2k+1}} )
		(|u_{2i}|^2 + |u_{2i+1}|^2) \tOmega^n \wedge \bOmega^n.
	\end{equation}
	For the fifth term of \eqref{eqn:C1less0}, we compute
	\begin{equation}
		\begin{split}
			\partial \partial_J u \wedge A &= \frac{n}{n-1} \partial \partial_J u
			\wedge (\frac{n \tOmega^{n-1} \wedge \Omega}{\Omega^n}\Omega^{n-1} -\tOmega^{n-1}) \\
			&= \frac{n}{n-1} (S_1(\partial \partial_J u)\Omega - \partial \partial_J u ) \wedge \tOmega^{n-1} \\
			&= n (\tOmega^n - \Omega_h \wedge \tOmega^{n-1}).
		\end{split}
	\end{equation}
	By compactness of $M$, there exists $\epsilon >0$ such that $\Omega_h \geq \epsilon \Omega$. Hence we obtain
	\begin{equation}
		\label{neq:C1-3}
		\begin{split}
			 \varphi' \frac{\partial \partial_J u \wedge A \wedge \bOmega^n}{\tOmega{\mathstrut}^n \wedge \bOmega^n} &= n\varphi'  -n\varphi' \frac{\Omega_h \wedge \tOmega^{n-1} \wedge \bOmega^n }{\tOmega{\mathstrut}^n \wedge \bOmega^n} \\
			&\leq n\varphi' -
			\epsilon  \varphi' \sum_{i=0}^{n-1} \frac{1}{\tOmega_{2i2i+1}}.
		\end{split}
	\end{equation}

	We  assume $\beta \gg 1$ otherwise we are finished.
	By \eqref{bt}, \eqref{eqn:C1-1-local}, \eqref{eqn:C1-2} and  \eqref{neq:C1-3}, the inequality \eqref{eqn:C1less0} becomes
	\begin{equation}
		\label{neq:C1-final}
		\begin{split}
			0 \leq &  - \frac{1}{\beta} \sum_{i=0}^{2n-1} (u_i(f)_{\overline{i}} + u_{\overline{i}}(f)_i) \\
			&+ \frac{(\varphi')^2 + \varphi''}{n -1}\sum_{i=0}^{n-1} (\sum_{k\neq i}\frac{1}{\tOmega_{2k2k+1}} )  (|u_{2i}|^2 + |u_{2i+1}|^2) \\
			&+ n\varphi' - (\epsilon\varphi' - C_1 \frac{\sum |u_j|}{\beta} - C_2\frac{\sum |u_{\overline{j}}|}{\beta})\sum_{i=0}^{n-1}\frac{1}{\tOmega_{2i2i+1}} - \varphi' \tilde u_t.
		\end{split}
	\end{equation}
	The first term is bounded from above. Now we take
	\begin{equation}
		\varphi(s) = \frac{\log(2s+C_0)}{2}.
	\end{equation}
	where $C_0$ is determined by $C^0$ estimate. Then  \eqref{neq:C1-final} becomes
	\begin{equation}
		\label{neq:C1-bound}
		C_3 \geq C_4\sum_{i=0}^{n-1} (\sum_{k\neq i}\frac{1}{\tOmega_{2k2k+1}})
		(|u_{2i}|^2 + |u_{2i+1}|^2) +
		C_5 \sum_{i=0}^{n-1}\frac{1}{\tOmega_{2i2i+1}}.
	\end{equation}
	Thus for any fixed $i$
	\begin{equation*}
		\tOmega_{2i2i+1} \geq \frac{C_5}{C_3} \geq C.
	\end{equation*}
	By equation \eqref{eqn:C1-main} we also have
	\begin{equation*}
		\frac{1}{\tOmega_{2i2i+1} } = e^{-u_t-f} \prod_{j\neq i}\tOmega_{2j2j+1} \geq \frac{C^{n-1}}{\sup_{M}e^{u_t + f}}, \ 0\le i\le n-1.
	\end{equation*}
	Then by \eqref{neq:C1-bound} we obtain  $\beta$ is uniformly bounded.
\end{proof}

\section{Bound on \texorpdfstring{$\partial \partial_J u$}{∂∂ju}}
\label{pC2}

\begin{theorem}
	\label{laplace-bound}
	Let $u$ be a solution  to \eqref{eqn:flow} on $M \times [0, T)$. Then there exists a constant $C$
	depending only on the fixed data $(I, J, K, g, \Omega, \Omega_h)$ and $f$ such that
	\begin{equation}
		\label{neq:pC2}
		\sup\limits_{M\times[0,T)}|\partial \partial_J u|_g \leq C.
	\end{equation}
\end{theorem}

\begin{proof}
	For simplicity denote
	\begin{equation*}
		\eta  =S_1(\p \p_J u).
	\end{equation*}
	Consider the function
	\begin{equation*}
		G = \log \eta - \varphi(\tilde u)
	\end{equation*}
	where  $\varphi$ is the same as before. 
	 For any $T_0\in (0,T)$, suppose $\max\limits_{M\times [0, T_0]}G=G(p_0, t_0)$ with $(p_0,t_0)\in M\times [0, T_0)$. We want to show $\eta(p_0,t_0)$ is uniformly bounded.
	We choose the normal coordinates around $p_0$. All the calculations are carried at $(p_0,t_0)$. We have
	\begin{equation*}
		\begin{split}
		0\le	\p_t G &= \frac{ \eta_t }{\eta} - \varphi' \tilde u_t,\\
			\p G &= \frac{\p \eta }{\eta} - \varphi' \p u = 0, \\
			\p_J G &= \frac{\p_J \eta}{\eta} - \varphi' \p_J u = 0, \\
			\p \p_J G &= \frac{\p \p_J \eta}{\eta} - ((\varphi')^2 + \varphi'')\p u \wedge \p_J u
			- \varphi' \p \p_J u.
		\end{split}
	\end{equation*}
	We further have
	\begin{equation}
		\label{neq:pC2-beginning}
		\begin{split}
			0 &\leq\mathcal{P}(G)=G_t -  \frac{\p \p_J G \wedge A \wedge \bOmega^n}{\tOmega{\mathstrut}^n \wedge \bOmega^n} \\
			&=\frac{ \eta_t }{\eta} - \varphi' \tilde u_t  - \frac{\p \p_J \eta \wedge A \wedge \bOmega^n}{\eta \tOmega{\mathstrut}^n \wedge \bOmega^n}
			+ ((\varphi')^2 + \varphi'')\frac{\p u \wedge \p_J u \wedge A \wedge \bOmega^n }{\tOmega{\mathstrut}^n \wedge \bOmega^n}
			+ \varphi' \frac{\p \p_J u \wedge A \wedge \bOmega^n}{\tOmega{\mathstrut}^n \wedge \bOmega^n}.
		\end{split}
	\end{equation}
	The last two terms were dealt with in the previous section. Since 
	\begin{equation*}
		\eta \Omega^n = n \p \p_J u \wedge \Omega^{n-1},
	\end{equation*}
	by taking $\p_t$ on both sides we have for $\eta_t$ in the first term
	\begin{equation}\label{etat}
		\eta_t = u_{t, p \overline{p}}.
	\end{equation}
	
	We now focus on $\p \p_J \eta$ in the third term of \eqref{neq:pC2-beginning}.
	By definition $\eta$ is real, and
	\begin{equation*}
		\eta \bOmega^n = n \bp \bp_J  u \wedge \bOmega^{n-1}.
	\end{equation*}
	Under the hyperK\"{a}hler condition $\dmath \Omega = 0$,  differentiating twice the above equation gives
	\begin{equation}
		\p \p_J \eta \wedge \bOmega^n
		= n  \p \p_J \bp \bp_J  u  \wedge \bOmega^{n-1}
		= n \bp \bp_J  \p \p_J  u  \wedge \bOmega^{n-1}.
	\end{equation}
	
	We know that (see \eqref{eq:ddju})
	\begin{equation*}
		\partial\partial_J u = (n-1)\Omega_h - S_1(\Omega_h)\Omega + S_1(\tOmega)\Omega - (n-1)\tOmega.
	\end{equation*}
	Thus
	\begin{equation}
		\bp \bp_J \p \p_J u = (n-1)\bp \bp_J \Omega_h -
		\bp \bp_J S_1(\Omega_h)\wedge \Omega +
		\bp \bp_J S_1(\tOmega)\wedge \Omega
		- (n-1)\bp \bp_J \tOmega,
	\end{equation}
	where we  used the hyperK\"{a}hler condition on $\Omega$. Now we have
	\begin{equation}
		\label{eqn:pC2-forth-diff}
		\begin{split}
			\p \p_J \eta \wedge A \wedge \bOmega^n
			&= n A \wedge  \bp \bp_J  \p \p_J u \wedge \bOmega^{n-1} \\
			&= n(n-1) A \wedge \bp \bp_J \Omega_h \wedge \bOmega^{n-1}
			- n\bp \bp_J S_1(\Omega_h) \wedge A \wedge \Omega \wedge \bOmega^{n-1}\\
			&\quad + n\bp \bp_J S_1(\tOmega)\wedge A \wedge \Omega \wedge \bOmega^{n-1}
			- n(n-1) A \wedge \bp \bp_J \tOmega \wedge \bOmega^{n-1}
		\end{split}
	\end{equation}
	Note that
	\begin{equation*}
		A \wedge \Omega = \frac{n}{n-1}\big(S_{n-1}(\tOmega)\Omega^n - \tOmega^{n-1}\big) \wedge \Omega
		=S_{n-1}(\tOmega)\Omega^n
	\end{equation*}
	and
	\begin{equation*}
		\bp \bp_J S_1(\tOmega) \wedge \Omega^n
		= n \bp \bp_J \tOmega \wedge \Omega^{n-1}.
	\end{equation*}
	The third term of \eqref{eqn:pC2-forth-diff} becomes
	\begin{equation*}
		\begin{split}
			n\bp \bp_J S_1(\tOmega)\wedge A \wedge \Omega \wedge \bOmega^{n-1}
			&=n \bp \bp_J S_1(\tOmega)\wedge (\Omega^n
			\cdot S_{n-1}(\tOmega)) \wedge \bOmega^{n-1} \\
			&= n^2 S_{n-1}(\tOmega) \bp \bp_J \tOmega
			\wedge \Omega^{n-1} \wedge \bOmega^{n-1}.
		\end{split}
	\end{equation*}
	The forth term is
	\begin{equation*}
		n(n-1)A \wedge \bp \bp_J \tOmega \wedge \bOmega^{n-1}
		= n^2S_{n-1}(\tOmega) \bp \bp_J \tOmega
		\wedge \Omega^{n-1} \wedge \bOmega^{n-1}
		- n^2\tOmega^{n-1} \wedge \bp \bp_J \tOmega \wedge \bOmega^{n-1}.
	\end{equation*}
	The first two terms of \eqref{eqn:pC2-forth-diff} are similar and we get
	\begin{equation*}
		\p \p_J \eta \wedge A \wedge \bOmega^n
		= n^2 \bp \bp_J \tOmega \wedge  \tOmega^{n-1} \wedge \bOmega^{n-1}
		- n^2 \bp \bp_J \Omega_h \wedge  \tOmega^{n-1} \wedge \bOmega^{n-1}
	\end{equation*}
	and
	\begin{equation}
		\label{eqn:pC2-target}
		\begin{split}
			\frac{\p \p_J \eta \wedge A \wedge \bOmega^n}{\eta \tOmega{\mathstrut}^n \wedge \bOmega^n}
			&= n^2\frac{
				\bp \bp_J \tOmega  \wedge \tOmega^{n-1}
				\wedge \bOmega^{n-1}}{
				\eta \tOmega{\mathstrut}^n \wedge \bOmega^n
			}
			- n^2\frac{
				\bp \bp_J \Omega_h  \wedge \tOmega^{n-1}
				\wedge \bOmega^{n-1}}{
				\eta \tOmega{\mathstrut}^n \wedge \bOmega^n
			} \\
			&= \frac{1}{\eta } \sum_{i=0}^{n-1} \sum_{p=0}^{2n-1}
			\frac{\tOmega_{2i2i+1,p\bar p}}{\tOmega_{2i2i+1}}
			- \frac{1}{\eta } \sum_{i=0}^{n-1} \sum_{p=0}^{2n-1}
			\frac{(\Omega_h)_{2i2i+1,p\bar p}}{\tOmega_{2i2i+1}}\\
			&\geq \frac{1}{\eta } \sum_{i=0}^{n-1} \sum_{p=0}^{2n-1}
			\frac{\tOmega_{2i2i+1,p\bar p}}{\tOmega_{2i2i+1}}
			- \frac{C_1}{\eta} \sum_{i=0}^{n-1} \frac{1}{\tOmega_{2i2i+1}}.
		\end{split}
	\end{equation}
	We now rewrite the right hand side of \eqref{eqn:pC2-target} using the equation
	\begin{equation}
		\label{eqn-log}
		\text{Pf} (\tOmega_{ij}) = e^{u_t + f} \text{Pf} (\Omega_{ij})
	\end{equation}
	where $\Omega^n = n! \text{Pf}(\Omega_{ij})dz^0 \wedge \cdots \wedge dz^{2n-1}$. Take logarithm of both sides
	\begin{equation}
		\label{eqn:pC2-log}
		\log \text{Pf} (\tOmega_{ij}) = u_t + f + \log \text{Pf} (\Omega_{ij}).
	\end{equation}
	Since $\bp \Omega = 0$, we have $\bp \text{Pf}(\Omega) = 0$. By taking $\bp$ of \eqref{eqn:pC2-log} and using $\text{Pf}(\tOmega_{ij})^2 = \det(\tOmega_{ij})$, we get
	\begin{equation}
		\frac{1}{2} \sum \tOmega^{ij} \tOmega_{ji,\bar p} =u_{t, \overline{p}} +  f_{\bar p}.
	\end{equation}
	By taking $\p$ of both sides we obtain
	\begin{equation}
		\label{eqn:pC2-2th-eq-diff}
		\frac{1}{2} \sum \tOmega^{ij} \tOmega_{ji,\bar p p}=
		\frac{1}{2} \sum \tOmega^{ik} \tOmega _{kl,p} \tOmega^{lj} \tOmega_{ji,\bar p} + f_{p \bar p} + u_{t, p\overline{p}} .
	\end{equation}
	In local coordinates, the left hand side of \eqref{eqn:pC2-2th-eq-diff} is
	\begin{equation}
		\label{eqn:pC2-eq-diff-left}
		\frac{1}{2} \sum \tOmega^{2i2i+1} \tOmega_{2i+12i,p \bar p}
		+ \frac{1}{2} \sum \tOmega^{2i+1 2i} \tOmega_{2i2i+1,p \bar p}
		= \sum \frac{\tOmega_{2i2i+1, p \bar p}}{\tOmega_{2i2i+1}}.
	\end{equation}
	It was proved in \cite{fu2023monge} that the first term of the right hand side of \eqref{eqn:pC2-2th-eq-diff} is nonnegative, i.e.
	\begin{equation}
		\label{neq:pC2-claim}
		\sum \tOmega^{ik} \tOmega _{kl,p} \tOmega^{lj} \tOmega_{ji,\bar p} \geq 0.
	\end{equation}
Hence we obtain
	\begin{equation}\label{eta1}
		\frac{\p \p_J \eta \wedge A \wedge \bOmega^n}{\eta \tOmega{\mathstrut}^n \wedge \bOmega^n}
		\geq \frac{1}{2 \eta  } \Delta_{I, g}f
		- \frac{C_1}{\eta} \sum_{i=0}^{n-1} \frac{1}{\tOmega_{2i2i+1}} + \frac{1}{\eta}u_{t,p\overline{p}}.
	\end{equation}
Inserting  \eqref{etat}, \eqref{eta1}, \eqref{eqn:C1-2} and \eqref{neq:C1-3} into  \eqref{neq:pC2-beginning}, we have
	\begin{equation}
		\label{neq:pC2-last}
		\begin{split}
			0 \leq -&\frac{1}{2 \eta } \Delta_{I, g}f
			+ \frac{(\varphi')^2 + \varphi''}{n-1}\sum_{i=0}^{n-1} (\sum_{k\neq i}\frac{1}{\tOmega_{2k2k+1}} )
			(|u_{2i}|^2 + |u_{2i+1}|^2) \\
			&+n\varphi' -
			\left( \epsilon \varphi'- \frac{C_1}{\eta} \right)
			\sum_{i=0}^{n-1} \frac{1}{\tOmega_{2i2i+1}} - \varphi' \tilde u_t.
		\end{split}
	\end{equation}
	Assuming $\eta \gg 1 $, we obtain from \eqref{neq:pC2-last}
	\begin{equation}
		C_2 \geq
		C_3 \sum_{i=0}^{n-1}\frac{1}{\tOmega_{2i2i+1}}.
	\end{equation}
	Hence all $\tOmega_{2i2i+1}$ are uniformly bounded.
	Since $\eta = S_1(\p \p_J u) = S_1(\tOmega) - S_1(\Omega_h) $, we can therefore obtain a uniform bound on $\eta$.

\end{proof}

\section{ Proof of Theorem \ref{thm:main}}
\label{C2a}
In \cite{TWWY15}, Tosatti-Wang-Weinkove-Yang derived $C^{2, \alpha}$ estimates for solutions of some nonlinear elliptic equations based on a bound on the Laplacian of the solution, which was improved and extended to parabolic equations by Chu \cite{Chu16}. Bedulli-Gentili-Venozzi \cite{bedulli2023parabolic} proved the $C^{2, \alpha}$ for the quaternionic complex Monge-Ampere equation. In this section we apply their techniques to derive the $C^{2, \alpha}$ estimates in our setting. Then the longtime existence and convergence follows.

We first need to rewrite equation \eqref{eqn:flow} in terms of real $(1, 1)$\nobreakdash-forms, which can be done by using the following relation
\begin{equation*}
	\frac{\Omega^n \wedge \bOmega^n}{(n!)^2} = \frac{\omega^{2n}}{(2n)!}.
\end{equation*}
And the equation is reformulated as
\begin{equation}
	\label{eqn:complex}
	\omega_u^{2n} = e^{2(u_t + f)} \omega^{2n},
\end{equation}
where $\omega$ and $\omega_u$ are induced by $\Omega$ and $\tOmega$ respectively. 

\begin{lemma}
	Let $u$ be a solution  to \eqref{eqn:flow} on $M \times [0, T)$ 
	and $\epsilon \in (0, T)$, then we have 
	\begin{equation}
		\label{neq: C2a}
		||\nabla^2 u||_{C^{\alpha}(M \times \interval[open right]{\epsilon}{T})} \leq C_{\epsilon, \alpha},
	\end{equation}
	where the constant $C_{\epsilon, \alpha} > 0$
	depending only on  $(I, J, K, g, \Omega, \Omega_h)$ , $f$, $\epsilon$ and $\alpha$.
\end{lemma}

\begin{proof}
The proof here follows from \cite{TWWY15}, \cite{Chu16} and \cite{Chu18}.  For any point $p \in M$, choose a local chart around $p$ that corresponds to the unit ball $B_1$ in $\mathbb{C}^{2n}$ with $I$\nobreakdash-holomorphic coordinates $(z^0, \dots, z^{2n-1})$. We have $\omega = \sqrt{-1}g_{i \bar j} dz^i \wedge d \bar z^j$ where $(g_{i \bar j}(x))$ is a positive definite $2n \times 2n$ hermitian matrix given by the metric at any point $x \in B_1$. We introduce the real coordinates by $z^i = x^i + \sqrt{-1}x^{2n + i}$ for $i = 0, \dots, 2n-1$.
 
 The complex structure $I$ corresponds to an endomorphism of the real tangent space which we still denote by $I$, written in matrix form 
 \begin{equation*}
 	I = 
 	\begin{pmatrix}
 		0 & -I_{2n} \\
 		I_{2n} & 0
 	\end{pmatrix},
 \end{equation*} 
 where $I_{2n}$ denotes the identity matrix. 
 
 For any $2n \times 2n$ hermitian matrix $H = A + \sqrt{-1}B$, the standard way to identify $H$ with a real symmetric matrix $\iota(H) \in  \operatorname{Sym}(4n)$ is defined as
 \begin{equation*}
 	\iota(H) =
 	\begin{pmatrix}
 		A & B \\
 		-B & A
 	\end{pmatrix}.
 \end{equation*}
 
 Let $Q_{(x, t)}(r)$ denote the domain $B_x(r)\times \interval[open left]{t-r^2}{t}$. We want to check the equation \eqref{eqn:complex} is of the following form as in \cite[p.~14]{Chu16}
 \begin{equation}
 	\label{eqn:model}
 	u_t(x, t) - F(S(x, t) + T(D_{\mathbb{R}}^2 u, x, t), x, t) = h(x, t) 
 \end{equation} 
 where $u$ is defined in $ Q_{(0,0)}(1)$ up to scaling and translation,  $D_\mathbb{R}^2 u$ is the real Hessian and the functions $F$, $S$ and $T$ are defined as the following.
 \begin{align*}
 &F : \operatorname{Sym}(4n) \times Q_{(0,0)}(1) \rightarrow \mathbb{R}, \quad F(N, x, t) := \frac{1}{2} \log \det(N);\\
 &S : Q_{(0,0)}(1) \rightarrow \operatorname{Sym}(4n), \quad S(x, t) := \iota(g_{i \bar j}(x));
 \end{align*}
 and 
 \begin{align*}
 	&T : \operatorname{Sym}(4n) \times Q_{(0,0)}(1) \rightarrow \operatorname{Sym}(4n), \\
 	&T(N, x, t) := \frac{1}{n-1}\Big( 
 	\frac{1}{8} \operatorname{tr}\big(\iota(g_{i \bar j}(x))^{-1}p(N)\big)\iota(g_{i \bar j}(x)) 
 	- G(N,x) 
 	\Big),
 \end{align*}
where
\begin{equation*}
	\begin{aligned}
		p(N) &:= \frac{1}{2}(N + {}^{t}INI), \\
		G(N, x) &:= \frac{1}{4}\big(p(N) + \iota({}^tJ(x))p(N)\iota(J(x))\big).
	\end{aligned}
\end{equation*}
Here we are using $J(x)$ as the matrix representation of  the complex structure $J$.
Observe that $p(D_\mathbb{R}^2 u) = 2\iota(D_{\mathbb{C}}^2 u)$, we have
\begin{equation*}
	G(D_\mathbb{R}^2 u, x) = \frac{1}{2} \Big(\iota(u_{i \bar j}) + 
	\iota(J)_i^{\bar k}\iota (D_{\mathbb{C}}^2 u)_{l \bar k}\iota(J)_{\bar j}^l \Big)(x)
	=\frac{1}{2} \iota \big( \operatorname{Re}(\p \p_J u(\cdot I, \cdot J))_{i \bar j} \big)(x).
\end{equation*}
Moreover, one can verify that
\begin{equation*}
	\operatorname{tr}\big(\iota(g_{i \bar j}(x))^{-1}p(D^2_{\mathbb{R}}u)\big) 
	= 4 \operatorname{tr}(g_{i \bar j}^{-1}(x) D_{\mathbb{C}}^2 u) = 4 \Delta_{I, g} u.
\end{equation*}
Notice that for a hermitian matrix $H$, $\det(\iota(H)) = \det(H)^2$, hence we get 
\begin{align*}
	& u_t(x, t) - F(S(x, t) + T(D_{\mathbb{R}}^2 u, x, t), x, t) \\
	=\;& \frac{1}{2} \log \det{\Big( \iota(g_{i \bar j}(x)) 
	+ \frac{1}{n-1}\big( (\frac{1}{2}\Delta_{I, g}u ) \iota(g_{i \bar j}(x)) -  
	\frac{1}{2} \iota \big( \operatorname{Re}(\p \p_J u(\cdot I, \cdot J))_{i \bar j} \big)(x)
	 \big) \Big) }\\
	 =\;& \log \det \Big( g_{i \bar j} (x) + \frac{1}{n-1}\big(S_1(\p \p_J u)g_{i \bar j} (x)  - 
	\frac{1}{2} \iota \big( \operatorname{Re}(\p \p_J u(\cdot I, \cdot J))_{i \bar j} \big)(x) \big) \Big)\\
	=\;& -2f(x) - \log \det(g_{i \bar j} (x)).
\end{align*}
Thus the equation \eqref{eqn:complex} is indeed of form \eqref{eqn:model}. 

It remains to verify that the functions $F$, $S$ and $T$ defined above satisfies all the assumptions $\mathbf{H1}$ to $\mathbf{H3}$ in \cite[p.~14]{Chu16}. From Theorem  \ref{laplace-bound} we have $tr_g g_u \leq C$, thus we get
\begin{equation*}
	C_0^{-1} I_{4n} \leq S(x, t) + T(D_{\mathbb{R}}^2 u, x, t) \leq C_0 I_{4n}.
\end{equation*}

Take the convex set $\mathcal{E}$ to be the set of matrices $N \in \operatorname{Sym}(4n)$ with
\begin{equation*}
	C_0^{-1} I_{4n} \leq N \leq C_0 I_{4n}.
\end{equation*}
It is straightforward that $\mathbf{H1}$, $\mathbf{H3}$ and $\mathbf{H2}.(1)$, $\mathbf{H2}.(2)$ hold. For $\mathbf{H2}.(3)$, we choose local coordinates such that $g(x) = Id$ and $J$ is block diagonal with only $J_{\overline{2i+1}}^{2i}$ and $J_{\overline{2i}}^{2i+1}$ non-zero, while $p(P)$ is diagonal with eigenvalues $\lambda_1, \lambda_1, \dots, \lambda_{2n}, \lambda_{2n} \geq 0$. 
Then one computes the eigenvalues of $T(P, x, t)$  are $ \frac{1}{2}\sum_{i \neq j}\lambda_i \geq 0$. Thus for $P \geq 0$ we have $T(P, x, t) \geq 0$, and let $K = 2(n-1)$, then $K^{-1}||P|| \leq ||T(P,x,t)|| \leq K||P||$.

Finally, to apply \cite[Theorem 5.1]{Chu16}, we need  overcome the lack of $C^0$ bound of $u$ using the same argument as in \cite[Lemma 6.1]{Chu18}. Specifically, we split into two cases $T < 1 $ and $T \geq 1$. If $T < 1$ then we have a $C^0$ bound on $u$ since by Lemma \ref{ut-bound} $\sup_{M\times \interval[open right]{0}{T} }|u_t| \leq C$. Hence Theorem 5.1 in \cite{Chu16} applies directly in this case.

If $T \geq 1$, for any $b \in (0, T-1)$, we consider
\begin{equation*}
	u_b(x, t) = u(x, t+b) - \inf\limits_{M \times \interval[open right]{b}{b+1} } u(x, t)
\end{equation*}
for all $t \in \interval[open right]{0}{1}$. By Lemma \ref{neq:tC0}, we have $\sup_{M \times [0, 1)} |u_b(x, t)| \leq C$.  Moreover, it is obvious that  $u_b$ also satisfies the equation,  thus we have a Laplacian bound on $u_b$. By Theorem 5.1 in \cite{Chu16} to $u_b$, for any $\epsilon \in (0, \frac{1}{2}) $, we have 
\begin{equation*}
	||\nabla^2 u||_{C^{\alpha}(M \times \interval[open right]{b + \epsilon}{b + 1})} 
	= 	||\nabla^2 u_b||_{C^{\alpha}(M \times \interval[open right]{\epsilon}{1})} 
	\leq C_{\epsilon, \alpha},
\end{equation*}
where $C_{\epsilon, \alpha}$ is a uniform constant depending only on the fixed data  $(I, J, K, g, \Omega, \Omega_h)$ , $f$, $\epsilon$ and $\alpha$. Since $b \in (0, T-1) $ is arbitrary, we obtain the estimate.
\end{proof}
\begin{proof}[Proof of Theorem \ref{thm:main}]
Once we have the $C^{2,\alpha}$ estimates, we obtain the longtime existence and the exponential convergence of $\tilde u$ similar as the argument in \cite{gill2011}. Let $\tilde u_{\infty}=\lim\limits_{t\rightarrow\infty} \tilde u(\cdot, t)$, then $\tilde u_{\infty}$ satisfies
\begin{equation*}
\begin{aligned}
	\big(\Omega_h + & \frac{1}{n-1}((\frac{1}{2}\Delta_{I, g}\tilde u_{\infty})\Omega - \partial \partial_J \tilde u_{\infty})\big)^n
	= e^{f + \tilde b} \Omega^n \\
	\Omega_h &+ \frac{1}{n-1}((\frac{1}{2}\Delta_{I, g}\tilde u_{\infty})\Omega - \partial \partial_J \tilde u_{\infty}) > 0,
\end{aligned}
\end{equation*}
where 
\begin{align*}
	\tilde b = \Big(\int_M \Omega^n\wedge\bOmega^n\Big)^{-1}\int_M\Big(
	\log{\frac{\big(\Omega_h + \frac{1}{n-1}((\frac{1}{2}\Delta_{I, g}\tilde u_{\infty})\Omega - \partial \partial_J \tilde u_{\infty})\big)^n}{\Omega^n}} -f\Big)\Omega^n\wedge\bOmega^n.
\end{align*}  
\end{proof}

\end{document}